\documentclass[10pt]{article}

\usepackage{textcase}
\usepackage[mathscr]{eucal}
\usepackage{amsmath,amssymb} 
\usepackage{graphicx}
\newcommand{\ben}{\begin{enumerate}}
\newcommand{\een}{\end{enumerate}}
\newcommand{\beqn}{\begin{equation}}
\newcommand{\eeqn}{\end{equation}}

\def\bpf{\noindent{\sc Proof.~}}
\def\epf{\, \hfill \qed}
\def\qed{\hfill \vbox{\hrule
 \hbox{\vrule\hbox to 5pt{\vbox to 8pt{\vfil}\hfil}\vrule}\hrule}}

\def\({\left(}
\def\){\right)}

\def\eps{\epsilon}
\def\p{\text{P}}

\newtheorem{theorem}{Theorem}[section]
\newtheorem{lemma}[theorem]{Lemma}

 \numberwithin{equation}{section}

   \title{On the Locality of the Pr\"ufer Code}
   
 \author{  Craig Lennon  
   \\  
   {\small Department of Mathematics} \\ {\small The Ohio State University} \\
{\small 231 W. 18$^\mathrm{th}$ Ave. } \\  [-0.8ex]{\small Columbus, OH 43210}  \\  [-0.8ex]
{\small \texttt{lennon.13@osu.edu}}  }

\begin{document}

 \maketitle

 \begin{abstract}
  The Pr\"ufer code is a bijection between trees on the vertex set  $[n]$ and strings on the set $[n]$ of length $n-2$ (Pr\"ufer strings of order $n$).  In this paper we examine the `locality' properties of the Pr\"ufer code, i.e. the effect of changing an element of the Pr\"ufer string on the structure of the corresponding tree.  
  Our measure for the distance between two trees $T,T^*$ is $\Delta(T,T^*)=n-1-\vert E(T)\cap E(T^*)\vert$.  We randomly mutate the $\mu$th element of the Pr\"ufer string of the tree $T$, changing it to the tree $T^*$, and we  asymptotically estimate the probability that this results in a change of $\ell$ edges, i.e. $P(\Delta=\ell\, \vert \, \mu).$ We find that $P(\Delta=\ell\, \vert \, \mu)$ is on the order of $ n^{-1/3+o(1)}$  for any integer $\ell>1,$ and  that $P(\Delta=1\, \vert \, \mu)=(1-\mu/n)^2+o(1).$
  This result implies that the probability of a `perfect' mutation in the Pr\"ufer code (one for which $\Delta(T,T^*)=1$) is $1/3.$
  
 \end{abstract}

 \section{Introduction}
 
 The Pr\"ufer code is a bijection between trees on the vertex set  $[n]:=\{ 1,\dots,n \}$ and strings on the set $[n]$ of length $n-2$ (which we will refer to as $P$-strings).  If we are given a tree $T$, we encode $T$ as a $P$-string as follows: 
   at step $i$ ($1\le i \le n-2$) of the encoding process the lowest number leaf is removed, and it's neighbor is recorded as $p_i$, the $i$th element of the $P$-string
   $$
P=(p_1,\dots,p_{n-2}), \quad p_i \in  [n], \quad (1 \le i \le n-2).
$$   
  We will describe a decoding algorithm in a moment.

First we observe that the Pr\"ufer code is one of many methods of representing trees as numeric strings,  \cite{DtoR}, \cite {GAreps}, \cite{Thom}.  A representation with the property that small changes in the representation lead to small changes in the represented object is said to have high locality, a desirable property when the representation is used in a genetic algorithm \cite{poorrep}, \cite {GAreps}.
 The distance between two numeric string tree representations is the number of elements in the string which differ, and the distance between two trees $T,T^*$ is measured by the number of edges in one tree which are not in the other:
 $$\Delta =\Delta^{(n)}=\Delta^{(n)} (T,T^*):= 
 n-1 - \vert E(T)\cap E(T^*) \vert, $$
 where $E(T)$ is the edge set of tree $T$.
 
   By a mutation in the $P$-string we mean the change of exactly one element of the $P$-string.  Thus we denote the set of all ordered pairs of P-strings differing in exactly one coordinate (the mutation space) by ${\cal M},$ and by ${\cal M}_\mu$ we mean the subset of the mutation space in which the P-strings differ in the $\mu$ th coordinate: 
$$
{\cal M} = \bigcup_{1 =\mu}^{ n-2} {\cal M}_\mu, \quad {\cal M}_\mu:= \left\{ (P,P^*) \, : \, p_i=p_i^* \text{ for } i \neq \mu, \text{ and } p_\mu \neq p_\mu^* \right\},
$$
where 
$$
P=(p_1,\dots,p_{n-2}), \quad P^*=(p_1^*,\dots,p_{n-2}^*),
$$
 so $\vert {\cal M}\vert=n^{n-2}(n-2)(n-1),$ and $\vert {\cal M}_\mu\vert=n^{n-2}(n-1)$.
We choose a pair $(P,P^*) \in {\cal M}$ uniformly at random, and the random variable $\Delta$ measures the distance between the trees corresponding to $(P,P^*)$.  Using $\p \(\{ \text{event} \} \vert \circ \)$ to denote conditional probability, we have
\begin{align*}
\p \( \Delta=\ell \) &= \sum_{\mu=1}^{n-2} \p\( \Delta=\ell\, \vert \,(P,P^*) \in {\cal M}_\mu \)\p \( (P,P^*) \in {\cal M}_\mu \) \\
&= \sum_{\mu=1}^{n-2} \p\( \Delta=\ell \, \vert \,(P,P^*) \in {\cal M}_\mu \)\frac{1}{n-2}.
\end{align*}
Hereafter we will represent the event $(P,P^*) \in {\cal M}_\mu$ by $\mu$, as in 
$$
  \p\( \{ \text{event} \}  \, \vert \, \mu \):= \p\( \{ \text{event} \}  \, \vert \,(P,P^*) \in {\cal M}_\mu \).
$$

   Computer assisted experiments conducted by Thompson (see \cite{Thom} page 195-196) for trees with a vertex size as large as $n=100$ led him to conjecture that:  
 
\beqn \label{Th1}
\lim_{n \to \infty} \p \( \Delta^{(n)}=1 \)=\frac{1}{3},
\eeqn
and that if $   \mu/n \to \alpha,$ then 
\beqn \label{Th2}
\lim_{n \to \infty} \p \( \Delta^{(n)} =1\, \big\vert \, \mu \) =(1-\alpha)^2.
\eeqn
      In a recent paper \cite{PS}, Paulden and Smith use combinatorial and numerical methods to develop conjectures about the exact value of $\p \( \Delta  =\ell \, \vert \, \mu \)$ for $\ell=1,2,$ and about the generic form that $\p \( \Delta =\ell\, \vert \, \mu \)$ would take for $\ell>2$.  These conjectures, if true, would prove \eqref{Th1}-\eqref{Th2}.  Unfortunately, the formulas representing the exact value of  $\p \( \Delta =\ell \, \vert \, \mu \)$ are complicated, even for $\ell=1,2$, and the proof of their correctness may be difficult.  
 In this paper we will show by a probabilistic method that \eqref{Th1}-\eqref{Th2} is  indeed correct, 
proving that 
\beqn \label{conj2}
\p \( \Delta^{(n)}=1\, \big\vert \, \mu \) =(1-\mu/n)^2+O\( n^{-1/3 }\ln^2n \),
\eeqn
and showing in the process that 
\beqn \label{conj3}
\p \( \Delta^{(n)}=\ell\, \big\vert \, \mu \) = O \(n^{-1/3 }\ln^2n \), \qquad (\ell>1).
\eeqn
 Of course  \eqref{conj2} implies \eqref{Th1}, because $\int_0^1 (1-\alpha)^2\, d \alpha=1/3.$
  In order to prove these results we will need to analyze the following $P$-string decoding algorithm, which we learned of from \cite{Cho}, \cite{PS}.

\subsection{A Decoding Algorithm}
     In the decoding algorithm, the $P$-string $P=(p_1,\dots,p_{n-2})$ is read from rear to front, so we begin the algorithm at step $n-2$ and count down to step $0$.  We begin a generic step $i$ with a tree $T_{i+1}$ which is a subgraph of the tree $T$ which was encoded as $P$.  This tree has vertex set $V_{i+1}$ of cardinality $n-i-1$ and edge set $E_{i+1}$ of cardinality $n- i-2$.  We will add to $T_{i+1}$ a vertex from $X_{i+1}:= [n] \setminus V_{i+1}$, and an edge, and the resulting tree $T_i$ will contain $T_{i+1}$ as a subgraph.  The vertex added at step $i$ of the decoding algorithm is the vertex which was removed at step $i+1$ of the encoding algorithm, and will be denoted by $y_i$. 
        A formal description of the decoding algorithm is given below.

  \begin{center}
   Decoding Algorithm  \\
    \rule[1mm]{6in}{.1mm}
  \end{center}

\noindent {\em Input:}  $P=(p_1,\dots,p_{n-2})$ and $X_{n-1}:=[n-1],$ $V_{n-1}=\{n\},\, E_{n-1}=\emptyset$, $p_{n-1}:=n$.\\
{\em Step $i$} ($1 \le i \le n-2$):  We begin with the set $X_{i+1}$ and a tree $T_{i+1}$ having vertex set $V_{i+1}$ and edge set $E_{i+1}$.  We examine entry $p_i$ of $P$.
\ben  \item If $p_i \in  X_{i+1},$ then set $y_i=p_i$ .
\item If $p_i \notin  X_{i+1},$ then let $y_i = \max X_{i+1} $ (the largest element of $ X_{i+1}$).
\een
In either case we add $y_i$ to the tree $T_{i+1},$ joining it by an edge to the vertex $p_{i+1}$ (which must already be a vertex of  $T_{i+1}$). So $X_{i }= X_{i+1}\setminus \{ y_i \}, V_{i }= V_{i+1}\cup \{ y_i \},$ and $E_{i}=E_{i+1} \cup \{ \,\{ y_i,p_{i+1}\} \,\}.$    

{ \em Step 0:}  We add $y_0,$ the only vertex in $X_1,$ and the edge $\{ y_0,p_1 \}$ to the tree $T_1$ to form the tree $T_0=T.$ 
\begin{center}
 \rule[1mm]{6in}{.1mm} 

\end{center}

In this algorithm, we do not need to know the values of $p_1,\dots,p_i$ until after step $i+1$.  We will take advantage of this by using the principle of deferred decisions.  With $\mu$ fixed, we will begin with $p_{\mu+1},\dots,p_{n-2}$ determined, but with $p_1,\dots, p_\mu,$ as yet undetermined.  We will then choose the values of the $p_i$ for $1 \le i \le \mu$ when the algorithm requires those values and no sooner.    

This will mean that the composition of the sets $ X_{i},V_{i},E_i$  will only be determined once we have conditioned on $p_i,\dots,p_{n-2}.$  When we compute the probability that $p_{i-1}$ is in a set ${\cal A}_i$ whose elements are determined by $p_j,\, j>i,$ (for example $X_{i}$ or $V_{i}$) we are implicitly using the law of total probability:
$$
\p \( p_{i-1}\in {\cal A}_i\, \vert \, \mu  \) = \sum_{P_i}  \p \( p_{i-1}\in {\cal A}_i \, \vert \, P_i\, ; \, \mu \)\p \(   P_i \, \vert \, \mu \),
$$ 
where the sum above is over all $P$-sub-strings $P_i=(p_i,\dots,p_{n-2}) $ of the appropriate length, and $\p \(    P_i \, \vert \, \mu \)$ is the probability of entries $i$ through $n-2$ of the $P$-string taking the values 
$(p_i,\dots,p_{n-2}).$  We will leave such conditioning as implicit when estimating probabilities  of the type $\p \( p_{i-1}\in {\cal A}_i\, \vert \, \mu  \).$

In the next section, we will use the principle of deferred decisions to easily find a lower bound for $\p \( \Delta =1\, \vert \, \mu \)$, and in later sections we will use similar techniques to establish asymptotically sharp upper bounds for $\p \( \Delta =1\, \vert \, \mu \)$, as well as for $\p \( \Delta =\ell \, \vert \, \mu \)$ ($\ell>1$).  The combination of these bounds will prove \eqref{conj2}-\eqref{conj3}.

\section{Lower Bounds } \label{upper}

For a fixed value of $\mu,$ we will construct a pair of strings from ${\cal M}_\mu $,  starting our construction with two partial strings 
$$
P_{\mu+1}=\(p_{\mu+1}, \dots, p_{n-2}\),\qquad P_{\mu+1}^*=\(p_{\mu+1}^*, \dots, p_{n-2}^*\), \quad p_{j} =p_{j}^*,
$$
where $p_j$ has been selected uniformly at random from $[n]$ for $\mu+1 \le j \le n-2.$  We have not yet chosen $p_j,p_j^*$ for $j\le \mu.$  We run the decoding algorithm from step $n-2$ down through step $\mu+1$, and at this point we have two trees $T_{\mu+1}=T_{\mu+1}^*$ as which $P_{\mu+1}=P_{\mu+1}^*$ have been partially decoded.  Of course we also have the sets $ V_{\mu+1}= V_{\mu+1}^*$ and $X_{\mu+1}=X_{\mu+1}^*$, where

$$
 V_{i }:= \{ j \, : \, j \text{ is a vertex of } T_{i } \}, \quad V_{i}^*:= \{ j \, : \, j \text{ is a vertex of } T_{i}^* \}, 
 $$
and $X_{i}=[n] \setminus V_{i}, X_{i}^*=[n] \setminus V_{i}^*$.  We let $E_i,\, E_i^*$ represent the edge sets of $T_{i},T_{i}^*.$

Now we choose $p_\mu $ and $p_\mu^* \neq  p_\mu,$ and execute step $\mu$ of the decoding algorithm.
There are two possibilities:
\ben
\item If both $p_\mu,p_\mu^* \in V_{\mu+1} \cup \{ \max  X_{\mu+1}     \} $,  then $y_i=y_i^* =\max   X_{\mu+1}   .$  We have added the same vertex and the same edge ($y_i$ and $\{ y_i,p_{\mu+1} \}$) to both $T_{\mu+1}$ and $T_{\mu+1}^*$.  We have $V_\mu=V_\mu^*$ and $E_\mu=E_\mu^*$.

\item One of $p_\mu,p_\mu^*$ is not an element of the set $ V_{\mu+1} \cup \{ \max  X_{\mu+1}   \}. $
\een

 We will denote the first of these two events by 
 \beqn \label{E}
 {\cal E}:= \{ \text{both $p_\mu,p_\mu^* \in V_{\mu+1} \cup \{ \max  X_{\mu+1}   \} $} \},
 \eeqn
     and we will show that on this event, $\Delta=1$ no matter what values of $p_j=p_j^*$ ($1 \le j \le \mu-1$) we choose to complete the strings $P,P^*$.  Thus
 $$
 {\cal E} \subseteq \{ \Delta =1  \} \Longrightarrow \p ({\cal E} \, \vert \, \mu) \le \p (\Delta =1 \, \vert \, \mu) .
 $$
Let us now prove the set containment shown in the previous line.

\bpf
Suppose that event ${\cal E}$ occurs, so that $V_\mu=V_\mu^*$ and $X_\mu=X_\mu^*$, and $T_\mu=T_\mu^*.$  Now choose $p_1,\dots, p_{\mu-1}$ uniformly at random from $[n],$ with $p_i^*=p_i$ for $1 \le i \le \mu-1.$  

At steps $\mu-1, \mu-2,\dots,0$ of the algorithm, we will, at every step, read the same entry $p_i=p_i^*$ from the strings $P,P^*.$  Because $X_\mu=X_\mu^*$  and $p_{\mu-1}=p_{\mu-1}^*$, the algorithm demands that we add to $T_\mu,T_\mu^*$ the same vertex $y_{\mu-1}=y_{\mu-1}^*.$  This in turn means that $X_{\mu-1}=X_{\mu-1}^*.$  In a similar fashion, for $0\le i \le \mu-2$ we have
$$
X_{i+1}=X_{i+1}^* \Longrightarrow y_{i}=y_{i}^*.
$$ 
 Thus at every step $i\le \mu$ of the algorithm we add the same vertex to  $V_{i+1},V_{i+1}^*.$  Furthermore, at every step we are adding the edge 
$\{ y_i,p_{i+1} \}$ to  $E_{i+1}$ and the edge $\{ y_i ,p_{i+1}^* \}$ to  $E_{i+1}^*$. 
 Since $p_i=p_i^*$ for $i \neq \mu$ and $p_\mu \neq p_\mu^*$, we add the same edge to $T_{i+1} $ 
 and $T_{i+1}^*$ at every step except at step $\mu-1$ at which we add  
 $\{y_{\mu-1},p_\mu  \}$ to $T_{\mu}$ and  $\{y_{\mu-1},p_\mu^*  \}$ ($\neq \{y_{\mu-1},p_\mu  \}$) to $T_{\mu}^*$.  Of course the same edge cannot be added to a tree twice, so at no point could we have added $\{y_{\mu-1},p_\mu^*  \}$ to $T$ or $\{y_{\mu-1},p_\mu   \}$ to $T^*$.  Thus $T$ and $T^*$ must have exactly $n-2$ edges in common, and  
$$
\Delta= \Delta^{(n)}(T,T^*):= n-1- \vert E(T)\cap E(T^*) \vert  =1.
$$
  \epf
  
  {\bf Note:}  We have proved that if $X_k =X_k^*$ for $k \le \mu$ then $X_j =X_j^*$ for all $j\le k$, that the same vertex is added at every step $j \le k,$ and that the same edge is added at every step $j\le \min \{k, \, \mu-2 \}.$  We will need this result later.\\
  
  Now we bound the conditional probability of event ${\cal E}$.
 
\begin{align}
 \p \( \Delta =1\, \vert \, \mu\) \ge \p \( {\cal E} \, \vert \, \mu \)&= \frac{n-\mu }{n } \cdot \frac{ n-\mu-1}{ n-1}  \nonumber
 \\   
 & = 1 - \frac{2 \mu}{n} +\frac{\mu^2}{n^2} +O \(n^{-1} \)
 .\nonumber
 \end{align}
 Thus we have
  \beqn \nonumber
 \p \(\Delta =1\, \vert \, \mu \) \ge (1-\mu/n)^2 +O \(n^{-1} \).
  \eeqn
  Of course $\p \( \{\Delta =\ell\}\cap {\cal E}   \, \vert \, \mu \)=0$ for $\ell>1$, so  in order to prove \eqref{conj2}-\eqref{conj3} it remains  to show that 
   \beqn \label{Delta=k}
   \p \( \{\Delta =\ell\}\cap {\cal E}^c  \, \vert \, \mu \)  =O \(n^{-1/3}\ln^2n\) , \quad (\ell \ge 1).\eeqn
       
   This endeavor will prove more complicated than the upper bounds, so we will need to establish some preliminary results and make some observations which will prove useful later.

  \section{Observations and Preliminary Results} \label{obs}
  
  Recall that after step $j$ of the decoding algorithm we have two sets $X_j,X_j^*$ of vertices which have not been placed in $T_j,T_j^*.$  For $j \ge \mu+1$, we know that $X_j=X_j^*,$ but we may have $X_j \neq X_j^*$ for $j\le \mu.$  So let us consider then the set ${\cal X}_j:=X_j\cup X_j^*$.  

  Our goal is to show that either ${\cal X}_j=X_j,$ or ${\cal X}_j$ consists of $X_j\cap X_j^*$ and of two additional vertices, one in $V_j \setminus V_j^*$ and one in $V_j^* \setminus V_j.$  This means  ${\cal X}_j$ has the following form:  
    \begin{align} \nonumber
{\cal X}_j:=&\{ x_1< \cdots <x_a< \min\{z_j,z_j^*  \}<x_{a+1}<\cdots<x_{a+b}< \\ 
& \qquad \qquad  \qquad \qquad  \qquad  \max \{z_j,z_j^*  \}< x_{a+b+c}< \cdots< x_{a+b+c}\}, \label{V} 
\end{align}
    where 
    $$ z_j \in V_j \setminus V_j^*, \quad z_j^* \in V_j^* \setminus V_j,\quad x_i \in X_j \cap X_j^*,\; \;(1 \le  i\le a+b+c),$$ and $a,b,c\ge 0$, with $a+b+c=j-1$.  We will consider a set ${\cal X}_j=X_j$ to also have the form shown above, but with $\{z_j,z_j^*  \}=\emptyset$ and $b(j)=c(j)=0,\, a(j)=j.$   Thus when showing that ${\cal X}_j $ is of the form \eqref{V}, our concern is to show that 1) there is at most one vertex $z_j \in V_j \setminus V_j^*,$ and 2) that there can be such a vertex if and only if there is exactly one vertex $z_j^* \in V_j^* \setminus V_j,$ so $\vert \{z_j,z_j^*  \}\vert$ is $0$ or $2$.

    For $j \ge \mu+1,$ the set ${\cal X}_j=X_j =X_j^*,$  and it is easy to see that ${\cal X}_\mu$ is of the form \eqref{V}.  Also, we showed in the previous section that if $X_k =X_k^*$ for $k\le \mu$ then $X_j =X_j^*$ for all $j\le k$.  Thus it is enough to show that if ${\cal X}_j$ ($j \le \mu$) is of the form \eqref{V} with $\{z_j,z_j^*  \} \neq \emptyset,$ then ${\cal X}_{j-1}$ is also of the form \eqref{V}. This will be shown in the process of examining what happens to a set ${\cal X}_j$ of the form \eqref{V} (with $\{z_j,z_j^*  \} \neq \emptyset$) at step $j-1$ of the decoding algorithm, an examination which will take most of this section.  In this examination we present  notation and develop results upon which our later probabilistic analysis will depend.  We begin by considering the parameters $a,b,c.$
 
 Of course, $$a=a(j),\quad b=b(j), \quad c=c(j),$$ depend on $j$, (and on $p_\mu^*$ and $p_i,$ $i\ge j$), but we will use  the letters $a,b,c$ when $j$ is clear.
  We let 
 $$
 A_j:= \{ x_1< \cdots <x_a \}, \quad B_j:= \{x_{a+1}<\cdots<x_{a+b}  \},
 $$
 and 
 $$
 C_j:= \{ x_{a+b+1}< \cdots< x_{a+b+c} \},
 $$
 so ${\cal X}_j=A_j\cup B_j \cup C_j \cup \{z_j,\,z_j^* \}.$
 
 Ultimately, we are interested not just in the set ${\cal X}_j$, but in the distance between two trees, i.e. $\Delta .$  We will find it useful to examine how this distance changes with each step of the decoding algorithm, so we define
 $$
 \Delta_j=\Delta_j^{(n)} \(T_j,T_j^*,T_{j+1},T_{j+1}^*\):= 1- \vert E_j\cap E_j^* \vert  + \vert E_{j+1} \cap E_{j+1}^* \vert, \quad (0\le j \le n-2),
 $$
 and observe that 
 \begin{align} \nonumber
  \Delta^{(n)} &=n-1-\vert E_0\cap E_0^* \vert  + \vert E_{n-1} \cap E_{n-1}^*\vert \\
  &=\Delta_0+\cdots +\Delta_{n-2} \label{deltaj}
  \end{align}
(recall that $ T_{n-1}$ is the single vertex $n$ and $T=T_0$).
     We add exactly one edge to each tree at each step of the algorithm, so the function $\Delta_j$ has a range in the set $\{ -1,0,1 \}.$  It is easy to check that $\Delta_{\mu }=1$ as long as $\min \{p_\mu , \, p_{\mu}^*\} \notin V_{\mu+1}\cup \{\max X_{\mu+1} \}$ (so on ${\cal E}^c$), and that $\Delta_{\mu -1}\ge 0$ (because $p_\mu \neq p_{\mu}^* $).  Further, if  $X_{j }=X_{j }^*$  and $j<\mu$, then we will add the same edge at every step $i<j$, so $\Delta_i=0$ for all $i < j.$

 Finally, we will need some notation to keep track of what neighbor a given vertex had when it was first added to the tree.  Thus for $v \in \{1,\dots,n-1 \}$ we denote by $ h(v)$ the neighbor of $v$ in $T_j$, where $j$ is the highest number such that $v$ is a vertex of $ T_j.$  Formally,
\beqn \label{h}
\text{for } v=y_{j }, \quad    h(v)=h_{ P }(v):= p_{j+1}  , \quad (P=(p_1,\dots,p_{n-2})).
\eeqn
  
    For example, if our string is $(4,3,2,2,7)$, then 
 $$  h(1)=4,\,h(2)=7,\,h(3)=2,\,h(4)=3,\,h(5)=2,\,h(6)=7.$$    
 Now we are prepared to examine the behavior of the parameters $a,b,c,$ and to make some crucial observations about the behavior of $\Delta_{j }$.  In the process we will show that if ${\cal X}_j$ is of the form \eqref{V} with $\{z_j,z_j^*  \} \neq \emptyset$ then ${\cal X}_{j-1}$ is of the same form (but possibly with $\{z_{j-1},z_{j-1}^*  \} = \emptyset$, meaning ${\cal X}_{j-1}=X_{j-1} $).  The observations below apply to all $1\le j \le \mu,$ except observations about the value of $\Delta_{j-1},$ which apply only to $j\le \mu-1.$  For $j \ge \mu $ we only need to remember that $\Delta_\mu=1$ on ${\cal E}^c$ and $\Delta_{\mu-1}\ge 0$.

\ben 

 \item \label{pinB}
 If $p_{j -1 } \in A_{j  }\cup B_{j  }\cup C_{j  }$, then $y_{j-1 }=y_{j-1 }^*=p_{j-1 } $, while $z_{j-1}=z_j,\,z_{j-1}^*=z_j^*,$ and $\Delta_{j-1}=0$ because we add the edge $\{  p_{j-1 },p_j\}$ to both of $T_j,T_j^*$.  
\ben 
 \item If $p_{j-1 } \in A_{j  }$ then  $a(j-1)=a(j)-1$, while $b(j-1)=b(j)$ and $c(j-1)=c(j)$.  
 \item If $p_{j-1 } \in B_{j  }$ then $b(j-1)=b(j)-1$ while $a(j-1)=a(j)$ and $c(j-1)=c(j)$. 
 \item If $p_{j-1 } \in C_{j  }$ then $c(j-1)=c(j)-1$ while $a(j-1)=a(j)$ and $b(j-1)=b(j).$

 \een
 Thus in every case, one of the parameters $a,b,c$ decreases by 1 while the others remain unchanged.

   \item \label{caseV} Suppose that $p_{j-1 } \in  {\cal V}_j:= V_j \cap V_j^* $.  Then 
 
 \ben
 \item  \label{case'} If $b(j)=c(j)=0$ then $y_{j-1 } =z_j^*$ and $y_{j-1 }^*=z_j$, so $X_{j-1 } = X_{j-1 }^*.$   While $\Delta_{j-1}$ could assume any of the values $-1,0,1,$ we have $\Delta_i=0$ for all $i<j-1.$  
 
  \item \label{s>0t=0}  First suppose that $ z_j<z_j^*  $ and  $b(j)> 0, c(j)=0$.  Then   
   $y_{j-1 }^*=x_{a+b } $ and $y_{j-1 }=z_j^*$, making $z_{j-1}^*=x_{a+b },$ $z_{j-1} =z_j.$  We have  $B_{j -1 }=B_{j  } \setminus \{ x_{a+b}  \}$, so  $a(j-1)=a(j),\,  b(j-1)=b(j)-1, \, c(j-1)=0 .$ Further, $\Delta_{j-1}=0$ if and only if the event 
  \beqn
  {\cal H}_{j-1}^*:= \{h_{P^*}(z_j^*)=p_j   \}
  \eeqn
    occurs, and otherwise $\Delta_{j-1}=1$.
    
    Similarly, if $ z_j>z_j^*$ and  $b(j)> 0, c(j)=0$, then $y_{j-1 }=x_{a+b } $ and $y_{j-1 }^*=z_j $ with $z_{j-1 }=x_{a+b } ,$ $z_{j-1}^* =z_j^*.$  The change in the values of $a,b,c$ are the same as in the case of $ z_j<z_j^*$.  We also have $\Delta_{j-1}=0$ if and only if the event 
  \beqn
  {\cal H}_{j-1} := \{h_{P}(z_j)=p_j^*   \}
  \eeqn
    occurs, and otherwise $\Delta_{j-1}=1$.  In summary, if $b(j)> 0, c(j)=0$ and $p_{j-1 } \in {\cal V}_j $, then $\Delta_{j-1}=1$ unless  ${\cal H}_{j-1}\cup{\cal H}_{j-1}^*$ occurs.

 \item If $b(j) \ge 0, c(j)>0$ and $p_{j-1} \in {\cal V}_j$ then $y_{j-1 }^*=y_{j-1 }=x_{a+b+c} \in C_j$, $z_{j-1} =z_j ,$ $z_{j-1}^* =z_j^*,$ and we have $a(j-1)=a(j),\, b(j-1)=b(j), \, c(j-1)=c(j)-1.$  Since we add the edge $\{ x_{a+b+c},p_{j} \}$ to both of $T_j,T_j^*$ we have $\Delta_{j-1}=0.$ 

    \een

 \item \label{maxz} Suppose that $p_{j-1 }=  \max \{ z_j,z_j^* \}  .$   
 
 \ben
 \item  \label{case} If  $b(j)=c(j)=0$  then the results are the same as in the case \ref{case'}.  
 
  \item \label{s>0t=0'}  If  $b(j)> 0, c(j)=0$  then the results are the same as in the case \ref{s>0t=0}.    
    
 \item \label{vj*}   Suppose $b(j) \ge 0, c(j)>0.$  If $z_j<z_j^*$ and $p_{j-1 } =z_j^*$ then $y_{j-1 }^*=x_{a+b+c} $ and $y_{j-1 }=z_j^*$, making $z_{j-1}^*=x_{a+b+c} ,$ $z_{j-1} =z_j .$  If $z_j>z_j^*$ and $p_{j-1 } =z_j $ then $y_{j-1 }=x_{a+b+c} $ and $y_{j-1 }^*=z_j$, making $z_{j-1}=x_{a+b+c}, $ $z_{j-1}^* =z_j^* .$  In both cases, $a(j-1)=a(j)$, but $B_{j -1 }=B_{j  }\cup C_{j  }\setminus \{ x_{a+b+c}  \}$, so $c(j-1)=0,\, b(j-1)=b(j)+c(j)-1.$ In this case we have $\Delta_{j-1} \ge 0$.

   \een

    \item \label{minz} The last remaining possibility is that  $p_{j-1 } = \min \{ z_j,z_j^* \}$.  
    \ben 
 
  \item \label{c=0} If   $c(j)=0$ then $y_{j-1 }=z_j^*$ and $y_{j-1 }^*=z_j $ so $X_{j-1 } = X_{j-1 }^*$. We have $\Delta_{j-1}\in \{-1,0,1\}$ and $\Delta_i=0$ for all $i<j-1.$

  \item \label{vj} If $ c(j)>0$ and $z_j<z_j^*$  then $y_{j-1 }=x_{a+b+c} $   and $y_{j-1 }^*=z_j$, making $z_{j-1}=x_{a+b+c} ,$ $z_{j-1}^*  =z_j^*.$  If $z_j>z_j^*$ then $y_{j-1 }^*=x_{a+b+c} $ and $y_{j-1 } =z_j^*$, making $z_{j-1}^*=x_{a+b+c} ,$ $z_{j-1}  =z_j.$  In both cases $a(j-1)=a(j)+b(j)$ because the set $A_{j -1 }= A_{j  }\cup B_{j  }$, and $B_{j -1 }= C_{j  }\setminus\{x_{a+b+c}   \}$, so $c(j-1)=0,\, b(j-1)=c(j)-1 .$  In this case we have $\Delta_{j-1}\ge 0$.

    \een
        
 \een
 
 We have shown that if ${\cal X}_{j }$ is of the form shown in \eqref{V} then  ${\cal X}_{j-1}$ will be of the same form.  Furthermore, if $\{z_{j }, z_{j }^*  \}\neq \emptyset ,$ then $\{z_{j-1}, z_{j-1}^*  \} = \emptyset$ (i.e. $X_{j-1} = X_{j-1}^*$)  can only occur if  $c(j)=0$, see cases \ref{case'}, \ref{case}, and \ref{c=0}.  In addition, we observe that $\vert {\cal V}_{j } \vert=n-j$ when
 $\{z_{j}, z_{j}^*  \} = \emptyset,$ and if $\{z_{j}, z_{j}^*  \} \neq \emptyset$ then $\vert {\cal V}_{j } \vert=n-j-1.$
   We have also seen that as $j$ {\em decreases}:  1) the parameter $c(j)$ never gets larger, and 2) the parameter $b(j)$ decreases by 1 if $p_{j-1} \in B_{j }$ and otherwise can only decrease if $p_{j-1} \in \{z_{j},z_{j}^* \}.$
    We end our analysis of the decoding algorithm with one last observation, which is that $\Delta_j=-1$ for at most one value of $j$, which is clear from an examination of cases \ref{case'}, \ref{case}, and \ref{c=0}, since only in these cases can $\Delta_j=-1$, and in every case $\Delta_i=0$ for all $i<j$.

   In light of the knowledge that $\Delta_j=-1$ at most once, that $\Delta_{\mu }=1$ on ${\cal E}^c,$ and of \eqref{deltaj},  we now see that (on ${\cal E}^c$) if there are $\ell+1$ indices $j_1,\dots j_{\ell+1}<\mu $ such that $\Delta_i=1$ (for all $i \in \{ j_1,\dots j_{\ell+1}  \}$), then $\Delta>\ell.$  Thus in order to show that $\Delta(T,T^*)>\ell$ it suffices to find $\ell+1$ such indices.  So we have reduced the `global' problem of bounding (from below) $\Delta =\Delta_0+\cdots +\Delta_{n-2}$ to the `local' problem of showing that it is likely (on ${\cal E}^c$) that for at least $\ell+1$ indices $i < \mu$ we have $\Delta_i=1$.  We will begin this process in the next section.

\section{Upper Bounds} \label{L}

We now begin the process of showing that  for any positive integer $\ell$, 
\beqn \label{Pell}
\p \( \{ \Delta=\ell \}\cap {\cal E}^c\, \vert \, \mu\) =O \(n^{-1/3}\ln^2n \).
\eeqn
 The event  ${\cal E}$ is the event that $p_{\mu },p_{\mu }^* \in V_{\mu+1}\cup \{ \max X_{\mu+1} \}  $, which is the event that ${\cal X}_\mu=X_\mu$ (equivalently $\{z_{\mu}, z_{\mu}^*  \}= \emptyset $).  So on  ${\cal E}^c$ we have $\{z_{\mu}, z_{\mu}^*  \}\neq \emptyset ,$ and ${\cal E}^c$ is the union of the following events: 
   \ben \item ${\cal E}_1:=\{  b(\mu)<  \delta_n  \} \cap \{ \{z_{\mu}, z_{\mu}^*  \}\neq \emptyset \}$,  \qquad $\delta_n=n^{1/3},$
   \item ${\cal E}_2:=\{   b(\mu)\ge \delta_n \}$,
   \een
   so 
   $$
   \p \( \{ \Delta=\ell \}\cap {\cal E}^c\, \vert \, \mu\) \le \p \(   {\cal E}_1\, \vert \, \mu\)+\p \( \{ \Delta=\ell \}\cap {\cal E}_2\, \vert \, \mu\).
   $$
Let us show now that 
\beqn \label{s<eps}
\p \( {\cal E}_1\, \vert \, \mu \) =O(\delta_n/n)     .
\eeqn
      \bpf
      Consider the sets $$
     {\cal X}_{\mu+1}= X_{\mu+1}\ =\{ x_1<\dots < x_{\mu+1} \},\qquad {\cal V}_{\mu+1}=V_{\mu+1} =[n]\setminus X_{\mu+1}.
      $$
      On ${\cal E}_1$ either:  1) $\max \{ p_{\mu },p_{\mu }^* \} \in {\cal V}_{\mu+1}$ and $\min\{ p_{\mu },p_{\mu }^* \} $ is one of the $\lfloor \delta_n \rfloor$ largest elements of ${\cal X}_{\mu+1},$ or 2) $  p_{\mu }  \in {\cal X}_{\mu+1}$ and $ p_{\mu }^*   $ is separated from $p_\mu$ by at most $\lfloor \delta_n \rfloor$ elements of ${\cal X}_{\mu+1}.$  So denote by ${\cal F} $ the event that $\max \{ p_{\mu },p_{\mu }^* \} \in {\cal V}_{\mu+1}$ and $\min\{ p_{\mu },p_{\mu }^* \} $ is one of the $\lfloor \delta_n \rfloor$ largest elements of ${\cal X}_{\mu+1}$. Then 
            $$
            {\cal F} \subseteq {\cal U}_1:= \{  \text{at least one of } p_{\mu },p_{\mu }^*   \text{ is one of the $\lfloor \delta_n \rfloor$ largest elements of ${\cal X}_{\mu+1}$} \}.
            $$
      Because $p_\mu$ is chosen uniformly at random from $[n]$ and $p_\mu^*$ is chosen uniformly at random from $[n]\setminus \{ p_\mu  \}$, a union bound gives us
      $$
      \p \(  {\cal F} \,\vert \, \mu \) \le  \p \(  {\cal U}_1 \,\vert \, \mu \) \le   \frac{\lfloor \delta_n \rfloor}{n} +\frac{\lfloor \delta_n \rfloor-1}{n-1}= O(\delta_n/n).
      $$
     On the event ${\cal E}_1\setminus {\cal F}$, we must have $p_\mu,p_\mu^* \in {\cal X}_{\mu+1}$ and there must be at most $\lfloor \delta_n \rfloor$ elements of ${\cal X}_{\mu+1}$ separating $p_\mu$ from $p_{\mu }^*.$  Thus we define
     \begin{align*}
     &{\cal U}_2:= \left\{ p_\mu =x_j \in X_{\mu+1}\, ; \, p_\mu^* \in {\cal Y}_j \right\},\\
     &{\cal Y}_j  :=\left\{  x_{\min \{ 1,j-\lfloor \delta_n \rfloor \}}, \dots,  x_{\max \{ \mu+1,j+\lfloor \delta_n \rfloor \}}  \right\} \setminus \{ x_j \}\subseteq X_{\mu+1}^*, \quad \vert {\cal Y}_j \vert \le 2 \lfloor \delta_n \rfloor
     \end{align*}
     and observe that ${\cal E}_1\setminus {\cal F}\subseteq {\cal U}_2.$  Then we have
      \begin{align*}    
\p \(  {\cal U}_2 \,\vert \, \mu \) &= \sum_{j=1}^{\mu+1} \p \( p_\mu^* \in {\cal Y}_j \,\vert \,p_\mu=x_j \, ;\, \mu \) \p \(p_\mu=x_j \in X_{\mu+1}\,\vert \,   \mu \) \\
& \le \sum_{j=1}^{\mu+1} \frac{2 \lfloor \delta_n \rfloor}{n-1} \frac{1}{n} =O(\delta_n/ n).      \end{align*}  
\qed

So 
we have proved \eqref{s<eps}, and from now on, we may assume that $b(\mu )=\vert B_\mu \vert$ is at least $\lceil \delta_n \rceil.$  Further, $B_\mu \subseteq X_j\setminus \{z_j  \},$ and $\vert X_\mu \vert = \mu,$ so we must have $\mu \ge \lceil \delta_n \rceil+1$ on the event ${\cal E}_2.$  So from here on we will also be restricting our attention to $\mu \ge \lceil \delta_n \rceil+1.$
 
 \subsection{The event ${\cal E}_2 $}
 
 In order to deal with ${\cal E}_2,$ we will  begin at step $\mu-1$, with $p_\mu^*,p_\mu,\dots,p_{n-2}$ already chosen, and we will begin choosing values for a number of positions $p_j=p_j^*$ ($j<\mu$) of our $P$-strings.    
 We will find that with high probability (whp) at some step $\tau=\tau(P,P^*) $ we have $c(\tau )=0,$ but $b(\tau )$ is on the order of $\delta_n$. 
So we will have at least $b(\tau ) $ values of $p_j$ ($j<\mu$) left to choose, and it is likely that for at least $\ell+1$ of those choices we will have $p_j \in {\cal V}_{j+1}.$   From case \ref{s>0t=0} of section \ref{obs}, we know that when this happens there are three possibilities:
 \ben 
 \item the event $ {\cal H}_{j }:=\{h_{P}(z_{j+1})=p_{j+1}^*   \}$ occurs,
 
 \item the event ${\cal H}_{j }^*:= \{h_{P^*}(z_{j+1}^*)=p_{j+1}   \} $ occurs, or 
 
 \item $\Delta_{j }=1$.  
\een
  The event $ {\cal H}_{j }\cup{\cal H}_{j }^*$ is unlikely to occur often, so (whp) we will have  $\Delta_{j }=1$ for at least $\ell+1$ values of $j<\mu $, which means that  $\Delta>\ell$ (whp).  \\
  
  To prove this, 
  let us define the random variable 
  $$\tau(z)=\tau(z)(P,P^*)= \max_{j   \le \mu } \{ j \, : \; c(j) \le z\} \qquad (\mu \ge \lceil \delta_n \rceil+1),   $$
and the events
\begin{align} 
 \label{T}
&{\cal S}:= \left\{b(\tau(0)) \ge  2^{-12} \delta_n  \right\}, \qquad \, \delta=\delta_n:= n^{1/3},   
\\
&{\cal T}_1 := \{ \tau(\delta) -\tau(0) \le 2 \beta_n  \}, \quad 
 {\cal T}_2 := \{  \tau(0) \le n- \beta_n \}, \quad \beta_n:=n^{2/3} \ln^2 n . \nonumber
\end{align}
We observe that for $u \le v$ we have $\tau(u)\le \tau(v)$ because $c(j)$ is a non-decreasing function of $j$ ($j\le \mu$).  Further, we note that if $\tau(z)< \mu,$ then $\vert C_{\tau(z)+1} \vert \ge z+1$, and because $C_j \subseteq X_j\setminus \{ z_j \},$ we have $\vert X_{\tau(z)+1} \vert \ge z+2.$  Since $\vert X_j \vert =j,$  it must be true that $\tau(z)\ge z+1,$ and in particular we have $\tau(\delta)\ge  \delta_n  +1,\, \tau(0)\ge 1.$  These bounds also hold if $\tau(\delta),\tau(0)=\mu.$  By a similar argument we can see that if $b(\tau(0))\ge  2^{-12} \delta_n$ (as on the event ${\cal S}$) then we must have $\tau(0)\ge 2^{-12} \delta_n+1.$
Finally, the following set containment holds for any sets ${\cal S},{\cal T}_1,{\cal T}_2$:
\beqn \label{Esub}
\{ \Delta=\ell \}\cap{\cal E}_2 \subseteq   {\cal T}_1^c \cup {\cal T}_2^c  \cup \( {\cal S}^c \cap {\cal T}_1\cap   {\cal E}_2 \) \cup  \(\{ \Delta=\ell \}\cap{\cal S}\cap{\cal T}_2 \)  .
\eeqn
In this section we will show first that 
\beqn \label{T1}
\p \( {\cal T}_1^c\, \vert \, \mu \)=O\( n^{-1} \),
\eeqn
second that
\beqn \label{T2}
\p \( {\cal T}_2^c\, \vert \, \mu \)=O(\beta_n/n),
\eeqn
 and finally that
\beqn \label{ST}
\p \( {\cal S}^c\cap{\cal T}_1\cap {\cal E}_2\, \vert \,\mu \)= O(\beta_n/n).
\eeqn
In section \ref{last} we will prove that 
\beqn \label{DS2}
  \p \( \{ \Delta=\ell\}\cap {\cal S}\cap{\cal T}_2 \, \vert\, \mu \) = O(\delta_n/n).
\eeqn
Combining results \eqref{T1}-\eqref{DS2} will prove, via \eqref{Esub}, that  
$$
\p\( \{ \Delta=\ell \}\cap {\cal E}_2 \, \vert \, \mu  \) = O(\beta_n/n)=O \(n^{-1/3} \ln^2 n  \).
$$

Since we are ultimately interested in the event $ \{ \Delta=\ell\}\cap {\cal S}\cap{\cal T}_2,$ which depends on $\tau(0)$, why must we concern ourselves with $\tau(\delta) $ and ${\cal T}_1$? 
To explain this, we must introduce the event
\begin{align} \label{defV}
&{\cal Z}_i:=\{  p_j \notin \{z_{j+1},z_{j+1}^* \}\text{ for $i \le j < \mu$}\}, \quad (1 \le i < \mu), \\
\nonumber 
&{\cal Z}_\delta:=\{  p_j \notin \{z_{j+1},z_{j+1}^* \}\text{ for $\tau(\delta) \le j < \mu$}\},\quad  \, {\cal Z}_0:=\{  p_j \notin \{z_{j+1},z_{j+1}^* \}\text{ for $\tau(0) \le j < \mu$}\}  .
\end{align}
For a fixed integer $i\ge 1,$ we know if the event ${\cal Z}_i$ occurred after examining $p_i,\dots,p_{n-2},p_\mu^*,$ while the events ${\cal Z}_\delta, \, {\cal Z}_0$ require knowledge of all $p_1,\dots ,p_{n-2},p_\mu^*.$  Of course if we condition on $\tau(0)$ or $\tau(\delta)$ then these last two events require knowledge of only $p_{\tau },\dots ,p_{n-2},p_\mu^*,$ for $\tau =\tau(0) ,\,\tau(\delta).$  Also, if $\tau(\delta)=\mu$ (respectively if $\tau(0)=\mu$) then the event ${\cal Z}_\delta$ (respectively ${\cal Z}_0$) trivially occurred.

 To see why we must consider $\tau(\delta),$ note that on the event 
$$\{p_j= \min \{z_{j+1},z_{j+1}^*  \} \text{ for } \tau(0)<j< \tau(\delta) \} \subseteq {\cal Z}_{0}^c$$
   we could have 
$$c(j+1)<<\delta_n \Longrightarrow b(j)=c(j+1)-1<<\delta_n, \quad c(j)=0,$$ see case \ref{vj} of section \ref{obs}.   This is a problem because we want $b(\tau(0)) $ to be at least on the order of $\delta_n.$  But if the event ${\cal Z}_{\delta}^c$ occurs, then 
 for some $j\ge \tau(\delta)$ either:  \\ $p_j= \min \{z_{j+1},z_{j+1}^*  \}$ and
$$c(j+1)\ge \delta_n+1 \Longrightarrow b(j)=c(j+1)-1 \ge \delta_n, \quad c(j)=0,$$ 
 (see case \ref{vj}), or $p_j= \max \{z_{j+1},z_{j+1}^*  \}$ and 
$$
c(j+1)\ge \delta_n+1 \Longrightarrow b(j)=b(j+1)+c(j+1)-1 \ge \delta_n, \quad c(j)=0,
$$
(see case \ref{vj*}).
 
Thus 
$$
{\cal Z}_\delta^c \subseteq {\cal S} \Longrightarrow 
{\cal S}^c \cap{\cal T}_1\cap {\cal E}_2 \subseteq \({\cal S}^c\cap {\cal Z}_0 \cap {\cal E}_2 \) \cup \( {\cal Z}_0^c\cap  {\cal Z}_\delta\cap{\cal T}_1 \),
$$
which means that
\beqn \label{STV}
\p\( {\cal S}^c \cap{\cal T}_1 \cap {\cal E}_2\, \vert \, \mu \)\le \p\( {\cal S}^c\cap {\cal Z}_0  \cap {\cal E}_2\, \vert \, \mu \)+\p\(   {\cal Z}_0^c\cap  {\cal Z}_\delta\cap{\cal T}_1 \, \vert \, \mu \).
\eeqn
In the process of proving \eqref{T1}, we will show that
\beqn \label{VT}
\p\(  {\cal Z}_0^c\cap{\cal Z}_\delta  \cap{\cal T}_1  \, \vert \, \mu \)=O(\beta_n/n),
\eeqn
and later in this section we will prove that 
\beqn \label{SV}
 \p\(  {\cal S}^c\cap {\cal Z}_0\cap {\cal E}_2\, \vert \, \mu \) =O\(n^{-1}\).
\eeqn
 The combination of \eqref{STV}-\eqref{SV}
 implies \eqref{ST}.
To conclude our remarks on the events ${\cal Z}_\delta, \, {\cal Z}_0,$ we note that an examination of their definitions shows that on ${\cal Z}_\delta$ (respectively on ${\cal Z}_0$) we cannot have reached $\tau(\delta)$ (resp. $\tau(0)$) by choosing $p_j \in \{ z_{j+1},z_{j+1}^* \}.$  Hence for $\tau(\delta)<\mu$ (resp. $\tau(0)<\mu$) we must have reached these points by choosing $p_j \in C_{j+1} \cup {\cal V}_{j+1},$ which in turn implies that the parameter $c(j) \ge c(j+1)-1$ for $j \ge \tau(\delta)$  (resp. $\tau(0)$).  On the other hand, on the set ${\cal Z}_\delta^c$ we have $\tau(\delta)=\tau(0).$

In the following proofs, we will occasionally show that $\p \({\cal B}\, \vert \,   \mu  \) \to 0$ by first showing that for some event ${\cal A}$ we have $\p \({\cal A}^c \, \vert \,   \mu  \) \to 0,$ and then showing that 
$$\p \({\cal B}\, \vert \, {\cal A}\, ; \, \mu  \):= \frac{\p \({\cal B}\cap {\cal A}\, \vert \, \mu  \)}{\p \(  {\cal A}\, \vert \,  \mu  \)} \to 0, \quad n \to \infty.$$
Obviously the result above proves that $  \p \({\cal B}\cap {\cal A}\, \vert \, \mu  \) \to 0$ as $n \to \infty.$  A conditional probability like the one above is only defined as long as $\p \(  {\cal A} \, \vert \, \mu\)>0,$ but of course if $\p \(  {\cal A} \, \vert \, \mu\)=0$ then because ${\cal B}\subset {\cal A} \cap {\cal A}^c$ we must have $\p \({\cal B}\, \vert \,   \mu  \) \to 0$ anyway.  Thus whenever we discuss conditional probabilities we will   assume (and not prove) that the event we condition on has positive probability.

 Let us begin proving the results we have discussed.
 \begin{lemma} 
\label{T1'} 
 Let   ${\cal T}_1= \{ \tau(\delta) -\tau(0) \le 2 \beta_n  \},$ and let ${\cal Z}_0,{\cal Z}_\delta$ be defined as in \eqref{defV}.  Then
 
 $$
  \p\(     {\cal T}_1^c  \, \vert \, \mu \) =O\(n^{-1}\),\qquad \quad \p\(  {\cal Z}_0^c\cap{\cal Z}_\delta  \cap{\cal T}_1  \, \vert \, \mu \)=O(\beta_n/n).
 $$
 \end{lemma}

\bpf  We will start with the second of the results above.  We will condition on the value of $\tau(\delta)$, and introduce notation for events conditioned on that value:  $$
 \p \( {\cal W}  \, \vert \, \tau \,;\,  \mu \) := \p \( {\cal W}  \, \vert  \tau=\tau(\delta) \,;\,  \mu \).
$$

With ${\cal Z}_i$ defined as in \eqref{defV}, we
  observe that ${\cal Z}_{ i } \subseteq {\cal Z}_{ i+1 }$.  If the set$\{z_{i+1}, z_{i+1}^*  \}$ is empty, then the (conditional) probability that $p_{i} \in \{z_{i+1}, z_{i+1}^*  \}$ is $0$, and if the set $\{z_{i+1}, z_{i+1}^*  \}$ is non-empty, and the (conditional) probability that $p_{i} \in \{z_{i+1}, z_{i+1}^*  \}$ is $2/n$.
  Thus we have
   \begin{align} \label{Zprod}
     \p \( {\cal Z}_i^c \cap {\cal Z}_{ i+1 }  \vert \, \tau \, ; \,   \mu\) \le 2/n, \quad ( 1 \le i <\mu-1).
      \end{align}
   To avoid having to condition also on the value of $\tau(0),$ we introduce ${\cal Z}_\phi$, where $\phi = \max \{ \tau(\delta)-2\lfloor  \beta_n   \rfloor , 0\}$ and note that with this definition, ${\cal Z}_\phi\cap {\cal T}_1 \subseteq {\cal Z}_0\cap {\cal T}_1$.  Also, a consideration of the definition of ${\cal Z}_i$ shows that on ${\cal Z}_\phi^c \cap {\cal Z}_\delta   $ we have $\tau(\delta)-\tau(0)\le 2\lfloor  \beta_n   \rfloor .$   
   
  From the law of total probability we have
 \begin{align} 
   \p \( {\cal Z}_\phi^c\cap  {\cal Z}_\delta     \, \vert \,   \mu \) &= \sum_{\tau=\lfloor \delta \rfloor+1}^{\mu}
       \p \( {\cal Z}_\phi^c\cap  {\cal Z}_\delta \, \vert \,   \tau \,;\,  \mu \)\p \(\tau=\tau(\delta) \,\vert \,  \mu \).
                  \label{V1}\end{align}
 Since $\tau -\phi \le 2\beta_n,$ 
 we obtain from \eqref{Zprod} the bound
 \begin{align}
   \p \( {\cal Z}_\phi^c \cap {\cal Z}_\delta  \, \vert \,   \tau \,;\,  \mu \) 
   &= \sum_{i=\phi}^{\tau -1}  \nonumber
   \p \( {\cal Z}_{i }^c\cap{\cal Z}_{i+1 } \, \vert \,   \tau \,;\,  \mu \) \\
 & \le 2\beta_n(2/n) = O(\beta_n/n). \label{VdV0}
 \end{align}
 This bound is independent of $\tau,$ so \eqref{VdV0}, combined with \eqref{V1} shows that 
\beqn \label{ZT}
 \p \( {\cal Z}_\phi^c\cap  {\cal Z}_\delta     \, \vert \,   \mu \)= O(\beta_n/n).
\eeqn

Because ${\cal Z}_\phi\cap {\cal T}_1 \subseteq {\cal Z}_0\cap {\cal T}_1,$   we have 
${\cal Z}_0^c\cap {\cal T}_1 \subseteq {\cal Z}_\phi^c\cap {\cal T}_1,$ so
 $$
  {\cal Z}_0^c\cap {\cal Z}_\delta  \cap {\cal T}_1 \subseteq {\cal Z}_\phi^c \cap {\cal Z}_\delta \cap  {\cal T}_1 
\subseteq {\cal Z}_\delta \cap {\cal Z}_\phi^c ,
$$
and \eqref{ZT} implies that $$
\p\(  {\cal Z}_0^c\cap{\cal Z}_\delta  \cap{\cal T}_1  \, \vert \, \mu \)=O(\beta_n/n).
$$

 Further, on the event ${\cal Z}_\delta^c$ we have $\tau(0)=\tau(\delta)$  and on the event ${\cal Z}_\delta \cap {\cal Z}_\phi^c$ we have $\tau(\delta)-\tau(0) \le 2\beta_n$, therefore 
 $${\cal Z}_\delta^c \subseteq {\cal T}_1, \; {\cal Z}_\delta \cap {\cal Z}_\phi^c \subseteq {\cal T}_1\Longrightarrow {\cal T}_1^c ={\cal T}_1^c \cap{\cal Z}_\phi.$$ 
 Thus 
 $$
 \p \({\cal T}_1^c \, \vert \, \mu \)= \p \({\cal T}_1^c \cap{\cal Z}_\phi \, \vert \, \mu \).
 $$

Now  $\{\tau(\delta)\le 2 \lfloor\beta_n \rfloor \} \subseteq {\cal T}_1,$ so when bounding the probability above we may restrict our attention to $\tau(\delta)> 2 \lfloor\beta_n \rfloor$.  Hence
\begin{align*}
\p \({\cal T}_1^c \cap{\cal Z}_\phi \, \vert \, \mu \) &= \sum_{\tau =2\lfloor  \beta_n \rfloor+1}^{n-2} \p \({\cal T}_1^c \cap{\cal Z}_\phi \, \vert \,  \tau \, ;\, \mu \) \p \( \tau(\delta)=\tau\, \vert \, \mu\) .
\end{align*}
To complete the proof of the lemma it is sufficient to show that
$$ 
\p \({\cal T}_1^c \cap{\cal Z}_\phi \, \vert \,  \tau\, ;\, \mu \)=O \(n^{-1}\).
$$
Toward this end we define  
 $$\eps_n:= \frac{1}{\ln n}, \quad \nu = \nu_n:=\lfloor  \eps_n \beta_n/\delta_n  \rfloor, \quad k=k_n:=\lfloor \delta_n/\eps_n \rfloor,$$
 observing that 
 $$
 k_n \nu_n \le  \beta_n, \qquad k_n >> \delta_n, \qquad k_n \nu_n^2 >>n \ln n.
 $$
 
 Then we 
consider the sub-string $(p_{\tau -2 \nu k},\dots,p_{\tau -1}),$ which can be divided into $2k$ segments of length $\nu$, leading us to introduce the notation
$$
 P(i):=(  p_{m(i)}, \dots, p_{m(i-1)-1}), \quad m(i):= \tau -i\nu, \quad (1 \le i \le 2k),
$$
and $${\cal D}_i:=\{  p_j \in {\cal V}_{j+1} \text{ for at least one } p_j \in  P(i)\} .
$$
The event ${\cal T}_1^c$ is the event that in steps $\tau-1$ through $\tau-\nu k$ we add fewer than $\delta_n$ elements of $C_{\tau(\delta)}$ as vertices of the pair of trees we are building.  Because every choice of a $p_j \in {\cal V}_{j+1}$ forces us to add a vertex from $C_{j+1}$, and because $k>>\delta_n$, we have 
  $$
 {\cal T}_1^c \subseteq \bigcup_{i=k+1}^{2k}{\cal D}_i^c   .
$$
So let us bound from above $ \p \({\cal D}_i^c\, \vert \, {\cal Z}_\phi\,;\,\tau \, ;\,\mu   \).$
 
On the event ${\cal Z}_\phi,$ we have 
$\vert {\cal V}_{j+1} \vert=n-(j+1)-1 $ for $\tau -2k\nu\le j \le \tau-1.$ 
Thus $$
 \p \( p_j \notin {\cal V}_{j+1} \, \vert \, {\cal Z}_\phi\,;\,\tau \, ;\,\mu   \)= 1-\frac{n-j-2}{n-2},
$$
and the events $p_j \notin {\cal V}_{j+1}$ are conditionally independent for $\tau -2k\nu\le j \le \tau-1.$
 Also for $m(i)\le j \le m(i-1)-1$ we have
\begin{align*}
\vert {\cal V}_{j+1} \vert=n-j-2 &\ge n-\big(\tau-(i-1) \nu-1 \big)-2 \\
& \ge n-\big(n-2-(i-1) \nu -1\big)-2 \\
&\ge(i-1) \nu.
\end{align*}
Thus we obtain the bound 
 
\begin{align}
 \p \({\cal D}_i^c\, \vert \, {\cal Z}_\phi\,;\,\tau \, ;\,\mu   \) &= \prod_{j=m(i)}^{ m(i-1)-1} \p \( p_j \notin {\cal V}_{j+1} \, \vert \, {\cal Z}_\phi\,;\,\tau \, ;\,\mu   \)\nonumber \\
  &= \prod_{j=m(i)}^{ m(i-1)-1}\( 1- \frac{n-j-2}{n-2 }\) \nonumber \\
  & \le  \( 1- \frac{(i-1)\nu}{n -2 } \)^\nu \le e^{-(i-1)\nu^2/(n-2)} .\label{Ac}
 \end{align}
   Hence 
\begin{align*}
\p \( {\cal T}_1^c\, \vert \,{\cal Z}_\phi\,;\, \tau \, ;\,\mu    \) &\le
\sum_{i=k+1}^{2k}\p \( {\cal D}_i^c \, \vert \,{\cal Z}_\phi\,;\, \tau \, ;\,\mu    \) \\
  &   \le k e^{-k\nu^2/(n-2)} =O\(n^{-1} \) ,
\end{align*}
and we find that
 $$
  \p \( {\cal T}_1^c\cap {\cal Z}_\phi\, \vert \,  \tau \, ;\,\mu    \)  \le \p \( {\cal T}_1^c \, \vert \, {\cal Z}_\phi\, ;\, \tau \, ;\,\mu    \)=O\(n^{-1}\) 
 .
 $$
 \qed 
  
 \begin{lemma} Let   ${\cal T}_2 = \{  \tau(0) \le n- \beta_n  \}$ and let ${\cal Z}_0,{\cal Z}_\delta$ be defined as in \eqref{defV}.  Then
 
 $$
  \p\(     {\cal T}_2^c  \, \vert \, \mu \) =O(\beta_n/n).
 $$
 \end{lemma}
\bpf Recall that by definition, $\tau(0)\le \mu,$ so the probability above is zero if $\mu \le n-\beta_n,$ and we may assume that $\mu \ge n-\beta_n.$  Now let us consider the set ${\cal Z}_{\rho},$ where $\rho=  \mu-\lfloor \beta_n \rfloor-1 ,$ and observe that on this event $c(j)\ge c(j+1)-1$ for $\rho  \le j <\mu.$  Thus 
$$
{\cal T}_2^c \subseteq {\cal Z}_{\rho}^c \cup \{c(\mu) < \lfloor\beta_n \rfloor+1\},
$$
and 
$$
 \p \({\cal T}_2^c   \, \vert \, \mu \)  \le  \p \({\cal Z}_{\rho}^c  \, \vert \, \mu \)  + \p \(\{c(\mu) < \lfloor\beta_n \rfloor+1\}  \, \vert \, \mu \) .
$$

We first observe that, by an argument similar to that in \eqref{VdV0}, we have
  \begin{align*}
    \p \({\cal Z}_{\rho}^c  \, \vert \, \mu \)   =O(\beta_n/n).
   \end{align*}
 Then we note that
$$
 \{c(\mu) < \lfloor \beta_n \rfloor+1\} \subseteq {\cal U}_1 \cup {\cal U}_2,
$$ 
where 
\begin{align*}
&{\cal U}_1 := \{ \max \{p_\mu, p_\mu^*\} \in V_{\mu+1} \}, \\
& {\cal U}_2 := \{ \max \{p_\mu, p_\mu^*\} \text{ is one of the $\lfloor \beta_n \rfloor+2$ largest elements of $X_{\mu+1}$} \}.
\end{align*}
So we have
$$
\p \( \{c(\mu) < \lfloor \beta_n \rfloor+1\}\, \vert \, \mu    \) \le 
\p \( {\cal U}_1\, \vert \, \mu     \) +
\p \( {\cal U}_2\, \vert \, \mu   \).
$$
Now for $\mu \ge n- \beta_n,$
\begin{align*}
\p \( {\cal U}_1\, \vert \, \mu  \)   \le      \frac{n-\mu-1}{n} +\frac{n-\mu-1}{n-1}  
  =O(\beta_n/n ),
\end{align*}
and
\begin{align*}
\p \( {\cal U}_2\, \vert \, \mu \)   \le   \frac{\lfloor \beta_n \rfloor+2}{n} +\frac{\lfloor \beta_n \rfloor+2}{n-1}     
 =O(\beta_n/n ).
\end{align*}
Thus 
$$
\p \( \{c(\mu) < \lfloor\beta_n \rfloor+1\}\, \vert \, \mu    \) =O(\beta_n/n ).
$$
\qed

\begin{lemma} Let ${\cal S}= \{b(\tau(0)) \ge  2^{-12} \delta_n  \},$  and let ${\cal Z}_0$ be defined as in \eqref{defV}.  Then
$$
\p \( {\cal S}^c\cap{\cal Z}_0\cap{\cal E}_2\, \vert \,\mu \) =O\(n^{-1}\).
$$

\end{lemma}

\bpf  Consider the event ${\cal Z}_0\cap \{ \tau(0)=\tau \}$.  On this event, if $\tau \le j$ then the only way we can have $b(j)<b(j+1)$ is if we choose $p_j \in B_{j+1}$, see section \ref{obs} case \ref{pinB}.  
 On the event  ${\cal E}_2\cap {\cal S}^c$ we have $b(\mu) \ge \delta_n$ but $b(\tau(0))< 2^{-12}\delta_n.$  Thus on the event ${\cal E}_2\cap {\cal S}^c \cap {\cal Z}_0  $  we must have chosen $p_j \in B_{j+1}$ more than $(1- 2^{-12})b(\mu)$ times over the range of indices $1\le j \le \mu-1.$  We will show that this   is unlikely to occur.

Toward this end, we will divide the substring $(p_1,\dots,p_{\mu-1})$ into segments again, this time letting 
$k(i)=\min\{ 0,\mu-in/12\}$, and for $i\ge 1$, we let 
$$
{\cal U}_i := B_\mu \cap \{x \notin \{p_{k(i )} , \dots,  p_{\mu-1} \}  \} , \quad u_i:= \vert {\cal U}_i\vert ,   \quad ({\cal U}_0:=B_\mu).
$$
So ${\cal U}_i $ (which depends on $p_{k(i )} , \dots,  p_{n-2},p_\mu^*$) is the set of elements of $B_\mu$ which have not been chosen as a $p_j$ for $j\ge k(i ) .$
We will show that with high probability $u_{i+1}\ge u_i/2$ for $0\le i \le 11,$ because if this happens for each such $i$ then we must have $u_{12} \ge 2^{-12}u_0.$  On the event ${\cal E}_2$, this implies the event ${\cal S}.$

Thus we have 
\begin{align} \label{EinJ}
{\cal E}_2\cap {\cal S}^c\cap {\cal Z}_0 &\subseteq 
  {\cal J}_{12}^c\cap {\cal E}_2,\qquad {\cal J}_i :=\bigcap_{j=0}^{i-1}  \{ u_{j+1} \ge u_j/2 \} , \quad (i\ge 1).
  \end{align}
 As \begin{align*}
 \p \({\cal J}_{12}^c\cap {\cal E}_2 \, \vert \, \mu \) &=\p \({\cal J}_1^c \cap {\cal E}_2\, \vert \, \mu \) + \sum_{i=2}^{12} \p \({\cal J}_i^c\cap {\cal J}_{i-1}\cap {\cal E}_2\, \vert \, \mu \)  \\
 & \le  \p \({\cal J}_1^c \, \vert \, {\cal E}_2  \, ;\, \mu \) + \sum_{i=2}^{12} \p \( {\cal J}_{i }^c  \, \vert \,{\cal J}_{i-1}\cap{\cal E}_2 \, ;\,  \mu\),
 \end{align*}
 it is enough to show that
\beqn \label{Ji}
\p \({\cal J}_1^c \, \vert \, {\cal E}_2  \, ;\, \mu \), \; \p \( {\cal J}_{i }^c  \, \vert \,{\cal J}_{i-1}\cap {\cal E}_2  \, ;\,  \mu\)=O\(n^{-1}\), \qquad (2\le i \le 12).
\eeqn
We will prove the result above for $\p \( {\cal J}_{i }^c  \, \vert \,{\cal J}_{i-1}  \cap {\cal E}_2\, ;\,  \mu\)$ -- the proof for $\p \({\cal J}_1^c \, \vert \, {\cal E}_2  \, ;\, \mu \)$ is similar.
Denote by  $\p \( \vert {\cal U}_i \vert =u \, \vert \,  {\cal J}_{i-1}\cap{\cal E}_2 \,;\, \mu\)$  the conditional probability that at the end of step $k(i ),$ the set ${\cal U}_i$ is a specific set of cardinality $u$ (${\cal U}_i=\{w_1, \dots, w_u  \}),$ and by $\p \( u_{i+1} < u /2\, \vert \,  {\cal J}_{i -1} \cap{\cal E}_2\,;\, \vert {\cal U}_i \vert =u\,;\, \mu \)$ the conditional probability that $u_{i+1} < u /2$ given the fixed set ${\cal U}_i$ (and given ${\cal J}_{i -1} \cap{\cal E}_2,\mu$).  Then   
 \begin{align*}
\p \( {\cal J}_{i }^c  \, \vert \,{\cal J}_{i-1} \cap {\cal E}_2 \, ;\,  \mu\) = \sum_{u=\lceil 2^{-i} \delta_n \rceil}^{n }\sum_{\vert {\cal U}_i \vert =u}  &\p \( u_{i+1} < u /2\, \vert \,  {\cal J}_{i -1} \cap{\cal E}_2\,;\, \vert {\cal U}_i \vert =u\,;\, \mu \)  \\
& \qquad \cdot \p \( \vert {\cal U}_i \vert =u \, \vert \,  {\cal J}_{i-1}\cap{\cal E}_2 \,;\, \mu\),
\end{align*}
where the outer sum above is over the cardinality of ${\cal U}_i$ and the inner sum is over all subsets of $[n ]$ of that cardinality.  The outer sum starts at $u=\lceil 2^{-i} \delta_n \rceil$ because conditioned on ${\cal J}_{i-1}\cap{\cal E}_2,$ we must have 
$$u \ge 2^{-i }b(\mu) \ge  2^{-i } \delta_n.$$
So we can prove  \eqref{Ji} by showing that 
\beqn \label{si}
\p \( u_{i+1} < u /2\, \vert \,  {\cal J}_{i -1} \cap{\cal E}_2\,;\, \vert {\cal U}_i \vert =u\,;\, \mu \)  =O\(n^{-1}\),
\eeqn
where the $O(\cdot)$ bound above is uniform over all sets ${\cal U}_i$ of cardinality at least $2^{-i} \delta_n$.  

The probability in \eqref{si} is equal to
$N({\cal U}_i)/n^{k(i)- k(i+1)}$, where
\ben \item  $
N({\cal U}_i)=$ the number of $P$-strings segments $(p_{k(i+1) },  \dots, p_{k(i)-1})$ such that we choose at least half of the elements of ${\cal U}_i$ as entries $p_j$ of our segment, and
\item  
  $n^{k(i)- k(i+1)}=$  the total number of $P$-strings segments $(p_{k(i+1) },  \dots, p_{k(i)-1})$. 
 \een
   Because we want to count $P$-strings segments, it is important that conditioning on the events ${\cal J}_{i-1}\cap{\cal E}_2$ and $\vert {\cal U}_i \vert =u$ requires knowledge of  $(p_{k(i ) },  \dots, p_{n-2}),p_\mu^*$ but not of the value of $p_j$ for $j\le k(i)-1$, and it is also important that for each $i$, $k(i)$ is a fixed number once we have conditioned on $\mu$.   
Before we begin counting, let us also introduce the notation 
$$(z)_j:=z(z-1) \cdots (z-j+1),  \quad  d=  \lfloor u/2 \rfloor ,$$
and note that for large enough $n$ we have $d  \ge 2^{-12}\delta_n.$
To find an upper bound for $N({\cal U}_i)$, we 
\ben \item choose $d$ out of $k:=k(i )- k(i+1)\le n/12$ positions,
\item choose $d$ distinct elements of  ${\cal U}_i $ for those positions, and 
\item then we choose any value of $p_j$ for the remaining 
 $k -d$ positions. 
 \een
   Thus, (for  $k-d \ge 0$) we have 
 \begin{align*}
 \frac{N({\cal U}_i)}{n^k}& \le {k \choose d } \frac{(u)_d \, n^{k -d}}{n^{k }} \\
 & 
= O\( \frac{ u^d  k^de^d}{d^d n^d}   \) \,\quad  \\& = O\(   2^de^d   12^{-d}\)=O\(2^{- \delta_n/2^{12}}\)=O\(n^{-1} \).
 \end{align*}
 For $k-d < 0,$ $N({\cal U}_i)=0$.  This proves \eqref{si}.  
 
   \qed

In this section we have shown  that 
$$
\p \(   \{ \Delta=\ell \} \cap {\cal E}_2 \, \vert \, \mu\) = \p \(   \{ \Delta=\ell \}\cap {\cal S}\cap {\cal T}_2 \, \vert \, \mu\)+O(\beta_n/n).
$$
  In the next section we will consider the event $ \{ \Delta=\ell \}\cap {\cal S}\cap {\cal T}_2.$

 \subsection{The event $ \{ \Delta=\ell \}\cap {\cal S}\cap {\cal T}_2$} \label{last}

Recall from case \ref{s>0t=0} of section \ref{obs} that if $b(j)>0,c(j)=0$ and we choose $p_{j-1 } \in   {\cal V}_j=  V_j \cap V_j^* $ then there are three possibilities:
 \ben 
 \item the event ${\cal H}_{j-1}:= \{h_{P}(z_j)=p_j^*   \}$ occured,
 
 \item the event $ {\cal H}_{j-1}^*:= \{h_{P^*}(z_j^*)=p_j   \}$ occured, or 
 
 \item $\Delta_{j-1}=1$.  
\een
On the event ${ \cal S},$  we have $b(\tau(0))\ge 2^{-12}\delta_n$ which implies that $\tau(0)\ge 2^{-12}\delta_n +1$ (see the discussion following  \eqref{T}).  Thus at step $\tau(0)$ we have at least $2^{-12}\delta_n$ values of $p_j$ ($j<\tau(0)$) left to choose, and we will show that it is likely that we will have $p_j \in {\cal V}_{j+1}$   at least $\ell+1$ times, and it is unlikely that $ {\cal H}_{j}^* , \, {\cal H}_{j}$ will occur for these $p_j$.  In this fashion we will show that  
\beqn \label{end}
\p \(\{ \Delta=\ell \}\cap {\cal S}\cap {\cal T}_2\, \vert \,  \mu \)=O\( \delta_n/n \).
\eeqn
To be more specific, we will let  
 $$\nu = \nu_n:=\lfloor   2^{-12}\delta_n/k  \rfloor, \quad k := \ell+1  , $$
and we will condition on the value of $\tau(0)$ ($\tau(0)=\tau$),  dividing the substring $(p_{\tau -k\nu}, \dots,p_{\tau -1})$,   into $k$ segments of length $\nu$, as we have done before. 
 We will find that this time we need to leave the first element of each segment as a buffer between adjacent  segments, so we use the notation
$$
P^- (i):= (p_{m(i)}, \dots,p_{m(i-1)-2} ) , \quad m(i):=\tau -i\nu,\quad (1 \le i \le k),
$$  
to denote the last $\nu-1$ elements of the $i$ th segment.  
On the event ${\cal T}_2=\{\tau(0)\le n -\beta_n    \}$ we have 
$$
  \vert {\cal V}_{j+1}\vert \ge n-j-2 \ge \beta_n-1, \quad ( 1\le j \le \tau -1 ).   
 $$
 Introducing the event \beqn \label{Z*}
 {\cal Z}_*:= \{ p_j \notin \{z_{j+1},z_{j+1}^*  \}   \text{ for }  \tau -k\nu \le j <\tau \} ,
 \eeqn
      we note that  
      $$
      \p\( p_j \in {\cal V}_{j+1}\, \vert \, {\cal Z}_* \cap {\cal T}_2\cap {\cal S}\, ; \, \tau \, ; \, \mu \) = \frac{n-j-2}{n-2}\ge \frac{\beta_n -1}{n-2},
      $$
      and that the events $ p_j \in {\cal V}_{j+1}$ are conditionally independent for $\tau -k\nu \le j <\tau.$
       We will show that the event ${\cal Z}_*^c$ is unlikely to occur conditioned on ${\cal T}_2\cap {\cal S}$.  Then we will find that, conditioned on ${\cal Z}_*$, it is likely that the event 
       $$
        { \cal C}:=\{\text{we choose at least one  $p_j \in {\cal V}_{j+1}$  in each segment $P^-(i)$} \}
       $$
       occurs.  At the same time we will prove a result which involves the buffer elements, i.e. for $\rho(i):=m(i-1)-2$ ($1 \le i \le k$) it is unlikely that the event ${\cal H}_{\rho(i)}\cup{\cal H}_{\rho(i)}^*$ will occur.  With all these results established, we will then be able to prove \eqref{end}.

\begin{lemma} \label{lemZ}

Conditioned on $\tau(0)=\tau$, let ${\cal Z}_*$ be defined as in \eqref{Z*}.  Then
$$
\p \({\cal Z}_*^c\,  \vert \, {\cal T}_2\cap {\cal S} \, ; \,\tau \, ; \, \mu \) =O(\delta_n/n).
$$
\end{lemma}

   \bpf   Let us begin by defining $$
{\cal Z}_*(i):= \{ p_j \notin \{z_{j+1},z_{j+1}^*  \} \text{ for } \,\tau - i \le j <\tau \} , \quad {\cal Z}_*(0):= \{ \{v_{\tau},v_{\tau }^*  \} \neq \emptyset \}.
$$

By the same argument as in \eqref{Zprod}, we have
$$\p \({\cal Z}_*(i)^c\cap {\cal Z}_*(i-1) \vert \,{\cal T}_2\cap {\cal S} \, ; \,  \tau\, ; \,  \mu \) \le 2/n , \qquad (1 \le i \le k \nu).$$ 
So
\begin{align*}
\p \({\cal Z}_*^c\,  \vert \, {\cal T}_2\cap {\cal S} \, ; \,\tau \, ; \, \mu \) &= \sum_{i=1}^{k\nu}
 \p \({\cal Z}_*(i)^c\cap {\cal Z}_*(i-1) \vert \,{\cal T}_2\cap {\cal S} \, ; \,  \tau\, ; \,  \mu \) \\
 & \le 2k\nu/n=O(\delta_n/n).
\end{align*}

\qed

Next, let 
  \begin{align*}
 { \cal C} := \bigcap_{i=1}^k {\cal C}_i, \qquad {\cal C}_i:= \{ p_j \in {\cal V}_{j+1} \text{ for at least one $p_j \in P^-(i)$}  \},
   \end{align*}  
and define
$$
H^{(\rho)}=H^{(\rho)}(P,P^*):=\sum_{i=1}^{k} {\bf I}_{{\cal H}_{\rho(i)}\cup{\cal H}_{\rho(i)}^*}, \quad \rho(i):=m(i-1)-2,
$$
   where $I_{\cal A}$ denotes the indicator of the event ${\cal A}.$  So $H^{(\rho)}$ counts the number of $i$ for which ${\cal H}_{\rho(i)}\cup{\cal H}_{\rho(i)}^*$ occurs.
   \begin{lemma} \label{lemC}
   Let ${\cal C},{\cal C}_i$ and $H^{(\rho)}$ be defined as above.  Then 
   $$
   \p \({ \cal C}^c\cup \left\{ H^{(\rho)} >0 \right\} \, \vert \,   {\cal Z}_* \cap{\cal T}_2\cap{\cal S} \, ; \, \tau \, ; \, \mu \)   =O\(n^{-1}  \).
   $$
   
   \end{lemma}
   
   \bpf 
 If we condition on ${\cal Z}_* \cap{\cal T}_2\cap{\cal S}$,  then for $\tau-k\nu \le j<\tau$ we have  
 \begin{align*}
 \vert {\cal V}_{j+1} \vert &= n-j-2 , \qquad \vert \{z_{j+1},z_{j+1}^* \} \vert =2,
   \end{align*}  
   and  the events $p_j \in {\cal V}_{j+1}$ are conditionally independent, with
   $$
   \p \( p_j \in {\cal V}_{j+1} \, \vert \,{\cal Z}_* \cap{\cal T}_2\cap{\cal S}\, ; \,\tau  \, ; \,  \mu \) = \frac{n-j-2}{n-2} \ge \frac{\beta_n-1}{n-2}.
   $$
  Thus, as in \eqref{Ac},  we obtain
 \begin{align*}
 \p \( {\cal C}_i^c\, \vert \,{\cal Z}_* \cap{\cal T}_2\cap{\cal S} \, ; \,\tau \, ; \,  \mu  \) &= \prod_{j=m(i)}^{\rho(i)} \(1-  \frac{n-j-3}{n-2}\) 
 \\ & \le \( 1- \frac{\beta_n-2}{n-2} \)^{\nu-1} =O\(  e^{- \beta_n \nu/n }\)=O\(n^{-1}\).
 \end{align*}
 
 Since $
 { \cal C}^c = \cup_{i=1}^k {\cal C}_i^c,  
 $
 we use a union bound to obtain 
 \begin{align*} \label{C1}
 \p ({ \cal C}^c\, \vert\,{\cal Z}_* \cap{\cal T}_2\cap{\cal S}\, ; \,\tau  \, ; \,  \mu ) \le   \sum_{i=1}^{k}  \p \({\cal C}_i^c \, \vert \,{\cal Z}_* \cap{\cal T}_2\cap{\cal S} \, ; \,\tau  \, ; \,  \mu \) = O\(n^{-1}\).
 \end{align*}
 Next we consider $H^{(\rho)}.$  
Conditioned on the event ${\cal Z}_* \cap{\cal T}_2\cap{\cal S}$, we have $\vert \{z_{j+1},z_{j+1}^* \} \vert =2$ (for $\tau-k\nu \le j<\tau$),
 so 
$$
 \p \(  {\cal H}_{\rho(i)}    \, \vert \, {\cal Z}_* \cap{\cal T}_2\cap{\cal S}  \,;\, \tau\,;\, \mu   \)= \left\{  
 \begin{aligned}
 &1/(n-2),\; \quad h_P(z_{\rho(i)+1}) \notin \{z_{\rho(i)+2},z_{\rho(i)+2}^* \}, \\
&0, \; \qquad \, \; \; \quad\quad \text{otherwise}.
 \end{aligned}
 \right.
$$
So we have
$$
 \p \(  {\cal H}_{\rho(i)}    \, \vert \, {\cal Z}_* \cap{\cal T}_2\cap{\cal S}  \,;\, \tau\,;\, \mu   \)=  \p \(  {\cal H}_{\rho(i)}^*    \, \vert \, {\cal Z}_* \cap{\cal T}_2\cap{\cal S}  \,;\, \tau\,;\, \mu   \)\le 
 \frac{1}{n-2},
 $$
and a union bound gives us 
$$
 \p \(  {\cal H}_{\rho(i)} \cup {\cal H}_{\rho(i)}^*  \, \vert \, {\cal Z}_* \cap{\cal T}_2\cap{\cal S}  \,;\, \tau\,;\, \mu   \) \le \frac{2}{n-2}.
$$
Hence
\begin{align*}
\p \( \left\{ H^{(\rho)} >0 \right\} \, \vert \,   {\cal Z}_* \cap{\cal T}_2\cap{\cal S} \, ; \, \tau \, ; \, \mu \)  
 &\le  \text{E} \left[H^{(\rho)}\, \vert \, {\cal Z}_* \cap{\cal T}_2\cap{\cal S} \,;\,\tau\,;\,\mu \right] \\
& =    \sum_{i=1}^{k} \p \(  {\cal H}_{\rho(i)}\cup{\cal H}_{\rho(i)}^*  \, \vert \, {\cal Z}_* \cap{\cal T}_2\cap{\cal S}\,;\, \tau\,;\, \mu   \) \\ &\le \frac{2k}{n-2}=O\( n^{-1}\).
\end{align*}
\qed

In order to complete our proof, we introduce the notation 
\beqn \label{G*}
{\cal G}={ \cal C}  \cap {\cal Z}_* \cap{\cal T}_2\cap{\cal S} \cap  \left\{ H^{(\rho)} =0 \right\} 
\eeqn
 and observe that 
$$
\{ \Delta=1 \}\cap {\cal S}\cap {\cal T}_2 \subseteq {\cal Z}_*^c\cup \(  \left[{ \cal C}^c\cup  \left\{ H^{(\rho)} >0 \right\}  \right] \cap {\cal Z}_* \cap{\cal T}_2\cap{\cal S} \) \cup\(\{ \Delta=1 \}\cap{\cal G} \).
$$
Lemmas \ref{lemZ} and \ref{lemC} imply
 \beqn \label{G}
   \p \( \{ \Delta=\ell \}\cap {\cal S}\cap {\cal T}_2\, \vert \,\tau  \, ; \, \mu\) =
  \p \( \{ \Delta=\ell \}\cap{\cal G} \, \vert \,\tau  \, ; \, \mu\)+O\( \delta_n/n \). 
\eeqn
so it remains only to prove the following lemma.
\begin{lemma}   
$$
\p \( \{ \Delta=\ell \}\cap{\cal G}  \, \vert \,  \tau  \, ; \, \mu\) =0.
$$
\end{lemma}
\bpf  On the event ${\cal G} $, we will choose at least one $p_j \in{\cal V}_{j+1}$ from each segment $P^-(i)$.  Thus   
  we can consider the (random) subset of indices
   \beqn \label{Gamma}
   \Gamma= \{ \gamma(1)< \dots <\gamma(k) \},   \eeqn for which $\gamma(i)$ is the largest element of $\{m(i), \dots ,   \rho(i)  \}$
      such that $ p_{\gamma(i)} \in {\cal V}_{\gamma(i)+1 }$.  This makes $p_{\gamma(i)}$ the last entry of the segment     such that  $p_j \in {\cal V}_{j+1  }.$

     We also define  
\beqn \label{CH}
 H^{(\gamma)}= H^{(\gamma)}(P,P^*):=  \sum_{i=1}^k {\bf I}_{{\cal H}_{\gamma(i)}\cup{\cal H}_{\gamma(i)}^*}, \quad     (  (P,P^*) \in {\cal G} ).
\eeqn
From the discussion at the beginning of this section, we can see that
$$
{\cal G} \cap \left\{ H^{(\gamma)} =0  \right\} \subseteq \{\Delta =\ell\}^c ,
$$
  which means that 
$$
\{\Delta =\ell\} \cap {\cal G}  \subseteq \left\{ H^{(\gamma)} >0  \right\}\cap  {\cal G}   \subseteq \left\{ H^{(\gamma)} >0  \right\} \cap \left\{H^{(\rho)}=0 \right\}.
$$
 To prove this lemma, it is enough to show that 
$$
\left\{ H^{(\gamma)} >0  \right\}\cap \left\{H^{(\rho)}=0 \right\}= \emptyset,
$$
which we can accomplish by proving that
\beqn \label{H<}
 H^{(\gamma)}(P,P^*) \le H^{(\rho)} (P,P^*) 
 \eeqn
   for all $(P,P^*) \in {\cal G}.$
 We begin by noting that, conditioned on ${\cal Z}_*,$ 
$$
p_j \in {\cal V}_{j+1}^c =  A_{j+1} \cup B_{j+1}\cup \{ z_{j+1},z_{j+1}^* \}   \Longrightarrow p_j \in A_{j+1} \cup B_{j+1},
$$
for $\tau-k \nu \le j<\tau$.  Thus if $  \gamma(i)=j<\rho(i)$, then 
$$
p_{j+1} \in A_{j+2} \cup B_{j+2}.
$$
Now, recall  that the elements of  $A_{j+2} \cup B_{j+2}$ have {\em not} appeared as any entry $p_i$ ($i \ge j+2$), but both $h_{P}(z_{j+1}),h_{P^*}(z_{j+1}^*)$ {\em have} appeared as some $p_i$ ($i \ge j+2$).  Thus 
$$
h_{P}(z_{j+1}),h_{P^*}(z_{j+1}^*) \notin A_{j+2} \cup B_{j+2} \Longrightarrow p_{j+1} \neq h_{P}(z_{j+1}),h_{P^*}(z_{j+1}^*).
$$ 
Consequently,
$$
\{ \gamma(i) <\rho(i)  \} \subseteq \( {\cal H}_{\gamma(i)}\cup{\cal H}_{\gamma(i)}^*  \)^c,
$$
which means that
$$
{\cal H}_{\gamma(i)}\cup{\cal H}_{\gamma(i)}^*\subseteq \( {\cal H}_{\rho(i)}\cup{\cal H}_{\rho(i)}^* \) \cap \{ \gamma(i) =\rho(i)  \} .
$$
So for every $i$ ($1 \le i \le k$),
$$
{\bf I}_{{\cal H}_{\gamma(i)}\cup{\cal H}_{\gamma(i)}^*} \le {\bf I}_{{\cal H}_{\rho(i)}\cup{\cal H}_{\rho(i)}^*} ,
$$
which proves \eqref{H<}.
\qed

\section{Conclusion}
In \cite{PS}, Paulden and Smith conjectured that $\p \(\Delta=\ell>1 \, \vert \, \mu \)$ was on the order of $n^{-1}$ (conjecture 3 on page 16).  We agree with this conjecture, even though we have only proved that $\p \(\Delta=\ell>1 \, \vert \, \mu \)$ is  on the order of $n^{-1/3+o(1)}.$  Our bound implies that  
$$
\lim_{n\to \infty}\p \(\Delta^{(n)}\ge  n^{ 1/3-o(1)}\, \Big\vert \, \mu \) = \frac{2}{3}.
$$
Thus, for large $n$, we should expect that a mutation in a $P$-string changes the structure of the tree by either one edge or by many edges, with little likelihood of anything in between occurring.

 \end{document}